\documentclass{amsart}
\usepackage{amssymb}

\usepackage{amscd}
\usepackage{thmdefs}

\input{tcilatex}
\theoremstyle{definition}
\theoremstyle{remark}
\numberwithin{equation}{section}

\input tcilatex

\begin{document}
\title[Generating Operator of $L(F_{N})$ ]{Moment Series and R-transform of the Generating Operator of $L(F_{N})$}
\author{Ilwoo Cho}
\address{Dep. of Math, Univ. of Iowa, Iowa City, IA, U. S. A}
\email{ilcho@math.uiowa.edu}
\keywords{Free Group Algebras, Amalgamated R-transforms, Amalgamated Moment Series,
Compatibility}
\maketitle

\begin{abstract}
In this paper, we will consider the free probability theory of free group
factor $L(F_{N}),$ where $F_{k}$ is the free group with $k$-generators. We
compute the moment series and the R-transform of the generating operator, $%
T=g_{1}+...+g_{N}+g_{1}^{-1}+..+g_{N}^{-1}.$
\end{abstract}

\strut

\strut

Voiculescu developed Free Probability Theory. Here, the classical concept of
Independence in Probability theory is replaced by a noncommutative analogue
called Freeness (See [9]). There are two approaches to study Free
Probability Theory. One of them is the original analytic approach of
Voiculescu and the other one is the combinatorial approach of Speicher and
Nica (See [1], [2] and [3]).\medskip Speicher defined the free cumulants
which are the main objects in Combinatorial approach of Free Probability
Theory. And he developed free probability theory by using Combinatorics and
Lattice theory on collections of noncrossing partitions (See [3]). Also,
Speicher considered the operator-valued free probability theory, which is
also defined and observed analytically by Voiculescu, when $\Bbb{C}$ is
replaced to an arbitrary algebra $B$ (See [1]). Nica defined R-transforms of
several random variables (See [2]). He defined these R-transforms as
multivariable formal series in noncommutative several indeterminants. To
observe the R-transform, the M\"{o}bius Inversion under the embedding of
lattices plays a key role (See [1],[3],[5],[12]).\strut In [12], we observed
the amalgamated R-transform calculus. Actually, amalgamated R-transforms are
defined originally by Voiculescu (See [10]) and are characterized
combinatorially by Speicher (See [1]). In [12], we defined amalgamated
R-transforms slightly differently from those defined in [1] and [10]. We
defined them as $B$-formal series and tried to characterize, like in [2] and
[3]. In [13], we observed the compatibility of a noncommutative probability
space and an amalgamated noncommutative probability space over an unital
algebra. In [14], we found the amalgamated moment series, the amalgamated
R-transform and the scalar-valued moment series and the scalar-valued
R-transform of the generating operator of $\Bbb{C}[F_{2}]*_{\Bbb{C}[F_{1}]}%
\Bbb{C}[F_{2}]$,

\strut

\begin{center}
$a+b+a^{-1}+b^{-1}+c+d+c^{-1}+d^{-1}$
\end{center}

\strut

where $<a,b>\,=F_{2}=\,<c,d>.$ The moment series and R-transforms
(operator-valued or scalar-valued) of the above generating operator is
determined by the recurrence relations. In this paper, by using one of the
recurrence relation found in [14], we will consider the moment series and
the R-transform of the ganarating operator of $L(F_{N}),$

\strut 

\begin{center}
$G=g_{1}+...+g_{N}+g_{1}^{-1}+...+g_{N}^{-1}\in L(F_{N}),$
\end{center}

\strut 

where $F_{N}=\,<g_{1},...,g_{N}>,$ for $N\in \Bbb{N}.$ 

\strut 

\strut 

\strut 

\section{Preliminaries}

\strut 

\strut 

Let $A$ be a von Neumann algebra and let $\tau :A\rightarrow \Bbb{C}$ be the
normalized faithful trace. Then we call the algebraic pair $(A,\tau ),$ the $%
W^{*}$-probability space and we call elements in $(A,\tau ),$ random
variables. Define the collection $\Theta _{s},$ consists of all formal
series without the constant terms in noncommutative indeterminants $%
z_{1},...,z_{s}$ ($s\in \Bbb{N}$). Then we can regard the moment series of
random variables and R-transforms of random variables as elements of $\Theta
_{s}.$ In fact, for any element \ $f\in \Theta _{s},$ there exists (some)
noncommutative probability space $(A,\tau )$ and random variables $%
x_{1},...,x_{s}\in (A,\tau )$ such that

\strut 

\begin{center}
$f(z_{1},...,z_{s})=R_{x_{1},...,x_{s}}(z_{1},....,z_{s}),$
\end{center}

\strut 

where $R_{x_{1},...,x_{s}}$ is the R-transform of $x_{1},...,x_{s},$ by Nica
and Speicher.

\strut 

\begin{definition}
Let $(A,\tau )$ be a $W^{*}$-probability space and let $a_{1},...,a_{s}\in
(A,\tau )$ be random variables ($s\in \Bbb{N}$). Define the moment series of 
$a_{1},...,a_{s}$ by

\strut 

$M_{a_{1},...,a_{w}}(z_{1},...,z_{s})=\sum_{n=1}^{\infty }\underset{%
i_{1},...,i_{n}\in \{1,...,s\}}{\sum }\tau \left(
a_{i_{1}}...a_{i_{n}}\right) z_{i_{1}}...z_{i_{n}}\in \Theta _{s}.$

\strut 

Define the R-transform of $a_{1},...,a_{s}$ by

\strut 

$R_{a_{1},...,a_{s}}(z_{1},...,z_{s})=\sum_{n=1}^{\infty }\underset{%
i_{1},...,i_{n}\in \{1,...,s\}}{\sum }k_{n}\left(
a_{i_{1}},...,a_{i_{n}}\right) \,z_{i_{1}}...z_{i_{n}}\in \Theta _{s}.$%
\strut 

\strut 

And we say that the $(i_{1},...,i_{n})$-th coefficient of $%
M_{a_{1},...,a_{s}}$ and that of $R_{a_{1},...,a_{s}}$ are the joint moment
of $a_{1},...,a_{s}$ and the joint cumulant of $a_{1},...,a_{s},$
respectively.
\end{definition}

\strut 

By Speicher and Nica, we have that

\strut 

\begin{proposition}
(See [1], [2] and [3]) Let $(A,\tau )$ be a noncommutative probability\
space and let $a_{1},...,a_{s}\in (A,\tau )$ be random variables. Then

\strut 

\begin{center}
$k_{n}\left( a_{i_{1}},...,a_{i_{n}}\right) =\underset{\pi \in NC(n)}{\sum }%
\tau _{\pi }\left( a_{i_{1}},...,a_{i_{n}}\right) \mu (\pi ,1_{n}),$
\end{center}

\strut 

where $NC(n)$ is the collection of all noncrossing partitions and $\mu $ is
the M\"{o}bius function in the incidence algebra and where $\tau _{\pi }$ is
the partition-dependent moment in the sense of Nica and Speicher (See [2]
and [3]), for all $(i_{1},...,i_{n})\in \{1,...,s\}^{n},$ $n\in \Bbb{N}.$
Equivalently,

\strut 

\begin{center}
$\tau \left( a_{i_{1}}...a_{i_{n}}\right) =\underset{\pi \in NC(n)}{\sum }%
k_{\pi }(a_{i_{1}},...,a_{i_{n}}),$
\end{center}

\strut 

where $k_{\pi }$ is the partition-dependent cumulant in the sense of Nica
and Speicher, for all $(i_{1},...,i_{n})\in \{1,...,s\}^{n},$ $n\in \Bbb{N}.$
$\square $
\end{proposition}

\strut 

The above combinatorial moment-cumulant relation is so-called the M\"{o}bius
inversion. The R-transforms play a key role to study the freeness and the
R-transform calculus is well-known (See [2] and [3]).\strut 

\strut 

Let $H$ be a group and let $L(H)$ be a group von Neumann algebra i.e

\strut 

\begin{center}
$L(H)=\overline{\Bbb{C}[H]}^{w}.$
\end{center}

\strut 

Precisely, we can regard $L(H)$ as a weak-closure of group algebra generated
by $H$ and hence

\strut 

\begin{center}
$L(H)=\overline{\{\underset{g\in H}{\sum }t_{g}g:g\in H\}}^{w}.$
\end{center}

\strut 

It is well known that $L(H)$ is a factor if and only if the given group $H$
is icc. Since our object $F_{N}$ is icc, the von Neumann group algebra $%
L(F_{N})$ is a factor. Now, define a trace $\tau :L(H)\rightarrow \Bbb{C}$
by 

\strut 

\begin{center}
$\tau \left( \underset{g\in H}{\sum }t_{g}g\right) =t_{e_{H}},$ \ for all \ $%
\underset{g\in H}{\sum }t_{g}g,$
\end{center}

\strut 

where $e_{H}$ is the identity of the group $H.$ Then $\left( L(H),\tau
\right) $ is the (group) $W^{*}$-tpobability space. Notice that $L(F_{N})$
is a II$_{1}$-factor under this trace $\tau .$\strut  Assume that the group $%
H$ has its generators $\{g_{j}\,:\,j\in I\}.$ We say that the operator

\strut 

\begin{center}
$G=\underset{j\in I}{\sum }g_{j}+\underset{j\in I}{\sum }g_{j}^{-1},$
\end{center}

\strut 

the generating operator of $L(H).$ 

\strut 

Rest of this paper, we will consider the moment series and the R-transform
of the generating operator $G$ of $L(F_{N}).$

\strut \strut 

\strut

\section{The Moment Series of the Generating Operator $G\in L(F_{N})$\strut }

\strut

\strut

In this chapter, we will consider free group II$_{1}$-factor, $A\overset{%
denote}{=}L(F_{N}),$ where $F_{k}$ is a free group with $k$-generators ($%
k\in \Bbb{N}$). i.e

\strut 

\begin{center}
$A=\{\underset{g\in F_{N}}{\sum }t_{g}g:t_{g}\in \Bbb{C}\}.$
\end{center}

\strut 

Recall that there is the canonical trace $\tau :A\rightarrow \Bbb{C}$
defined by

\strut

\begin{center}
$\tau \left( \underset{g\in F_{N}}{\sum }t_{g}g\right) =t_{e},$
\end{center}

\strut

where $e\in F_{N}$ is the identity of $F_{N}$ and hence $e\in L(F_{N})$ is
the unity $1_{L(F_{N})}.$ The algebraic pair $\left( L(F_{N}),\tau \right) $
is a $W^{*}$-probability space. Let $G$ be the generating operator of $%
L(F_{N}).$ i.e

\strut 

\begin{center}
$G=g_{1}+...+g_{N}+g_{1}^{-1}+...+g_{N}^{-1},$
\end{center}

\strut 

where \ $F_{N}=\,<g_{1},...,g_{N}>.$  It is known that if we denote $X_{n}=%
\underset{\left| w\right| =n}{\sum }w\in A$, for all $n\in \Bbb{N},$ then

\strut 

\ \ \ \ \ \ \ \ \ \ (1.1)

\begin{center}
$X_{1}X_{1}=X_{2}+2N\cdot e$ \ \ \ \ ($n=1$)
\end{center}

and

\ \ \ \ \ \ \ \ \ \ (1.2)

\begin{center}
$X_{1}X_{n}=X_{n+1}+(2N-1)X_{n-1}$ $\ \ \ \ (n\geq 2)$
\end{center}

\strut 

(See [36]).

In our case, we can regard our generating operator $G$ as $X_{1}$ in $A.$ 

By using the relation (1.1) and (1.2), we can express $G^{n}$ in terms of $%
X_{k}$'s ; For example, $G=X_{1},$

\strut

$G^{2}=X_{1}X_{1}=X_{2}+2N\cdot e,$

\strut

$G^{3}=X_{1}\cdot X_{1}^{2}=X_{1}\left( X_{2}+(2N)e\right)
=X_{1}X_{2}+(2N)X_{1}$

$\ \ \ \ =X_{3}+(2N-1)X_{1}+(2N)X_{1}=X_{3}+\left( (2N-1)+2N\right) X_{1},$

\strut

continuing

\strut

$G^{4}=X_{4}+\left( (2N-1)+(2N-1)+2N\right) X_{2}+(2N)\left(
(2N-1)+(2N)\right) e,$

$G^{5}=X_{5}+\left( (2N-1)+(2N-1)+(2N-1)+2N\right) X_{3}$

$\ \ \ \ \ \ \ \ \ +\left( (2N-1)\left( (2N-1)+(2N-1)+(2N)\right)
+(2N)\left( (2N-1)+(2N)\right) \right) X_{1},$

$G^{6}=X_{6}+\left( (2N-1)+(2N-1)+(2N-1)+(2N-1)+2N\right) X_{4}$

$\,\,\,\ \ \ \ \ \ \ \ \ +\{(2N-1)\left( (2N-1)+(2N-1)+(2N-1)+(2N)\right) $

$\ \ \ \ \ \ \ \ \ \ \ \ \ \ \ \ \ \ \ \ \ \ \ \ \ \ \ \ \ \ \ \
+(2N-1)\left( (2N-1)+(2N-1)+(2N)\right) $

$\ \ \ \ \ \ \ \ \ \ \ \ \ \ \ \ \ \ \ \ \ \ \ \ \ \ \ \ \ \ \ \ \ \ \ \ \ \
\ \ \ \ \ \ \ \ \ \ \ \ \ \ +(2N-1)((2N-1)+(2N))\}X_{2}$

$\ \ \ \ \ \ \ \ \ \ +(2N)\left( (2N-1)\left( (2N-1)+(2N-1)+(2N))\right)
+(2N)((2N-1)+(2N))\right) e,$

etc.

\strut

So, we can find a recurrence relation to get $G^{n}$ ($n\in \Bbb{N}$) with
respect to $X_{k}$'s ($k\leq n$). Inductively, we have that $G^{2k-1}$ and $%
G^{2k}$ have their representations in terms of $X_{j}$'s as follows ;

\strut

\begin{center}
$%
G^{2k-1}=X_{1}^{2k-1}=X_{2k-1}+q_{2k-3}^{2k-1}X_{2k-3}+q_{2k-5}^{2k-1}X_{2k-5}+...+q_{3}^{2k-1}X_{3}+q_{1}^{2k-1}X_{1}
$
\end{center}

\strut \strut

and

\begin{center}
$%
G^{2k}=X_{1}^{2k}=X_{2k}+p_{2k-2}^{2k}X_{2k-2}+p_{2k-4}^{2k}X_{2k-4}+...+p_{2}^{2k}X_{2}+p_{0}^{2k}e,
$
\end{center}

\strut

where $k\geq 2.$ Also, we have the following recurrence relation ;

\strut \strut 

\begin{proposition}
Let's fix $k\in \Bbb{N}\,\setminus \,\{1\}.$ Let $q_{i}^{2k-1}$ and $%
p_{j}^{2k}$ ($i=1,3,5,...,2k-1,....$ and $j=0,2,4,...,2k,...$) be given as
before. If $p_{0}^{2}=2N$ and $q_{1}^{3}=(2N-1)+(2N)^{2},$ then we have the
following recurrence relations ;

\strut 

(1) Let

\begin{center}
$%
G^{2k-1}=X_{2k-1}+q_{2k-3}^{2k-1}X_{2k-3}+...+q_{3}^{2k-1}X_{3}+q_{1}^{2k-1}X_{1}.
$
\end{center}

Then

\strut 

$G^{2k}=X_{2k}+\left( (2N-1)+q_{2k-3}^{2k-1}\right) X_{2k-2}+\left(
(2N-1)q_{2k-3}^{2k-1}+q_{2k-5}^{2k-1}\right) X_{2k-4}$

$\ \ \ \ \ \ \ \ \ \ \ \ \ \ \ \ \ +\left(
(2N-1)q_{2k-5}^{2k-1}+q_{2k-7}^{2k-1}\right) X_{2k-6}+$

$\ \ \ \ \ \ \ \ \ \ \ \ \ \ \ \ \ +...+\left(
(2N-1)q_{3}^{2k-1}+q_{1}^{2k-1}\right) X_{2}+(2N)q_{1}^{2k-1}e.$

\strut i.e,

\strut 

$\ \ \ \ \ \ \ p_{2k-2}^{2k}=(2N-1)+q_{2k-3}^{2k-1},$ $\ $

$\ \ \ \ \ \ \ p_{2k-4}^{2k}=(2N-1)q_{2k-3}^{2k-1}+q_{2k-5}^{2k-1},$

.....$...,$

\ \ \ \ \ \ \  $p_{2}^{2k}=(2N-1)q_{3}^{2k-1}+q_{1}^{2k-1}$ 

and 

$\ \ \ \ \ \ \ \ p_{0}^{2k}=(2N)q_{1}^{2k-1}.$

\strut 

(2) Let

\begin{center}
$G^{2k}=X_{2k}+p_{2k-2}^{2k}X_{2k-2}+...+p_{2}^{2k}X_{2}+p_{0}^{2k}e.$
\end{center}

Then

\strut 

$G^{2k+1}=X_{2k+1}+\left( (2N-1)+p_{2k-2}^{2k}\right) X_{2k-1}+\left(
(2N-1)p_{2k-2}^{2k}+p_{2k-4}^{2k}\right) X_{2k-3}$

$\ \ \ \ \ \ \ \ \ \ \ \ \ \ \ \ \ \ \ \ \ +\left(
(2N-1)p_{2k-4}^{2k}+p_{2k-6}^{2k}\right) X_{2k-5}+$

$\ \ \ \ \ \ \ \ \ \ \ \ \ \ \ \ \ \ \ \ \ +...+\left(
(2N-1)p_{4}^{2k}+p_{2}^{2k}\right) X_{3}+\left(
(2N-1)p_{2}^{2k}+p_{0}^{2k}\right) X_{1}.$

i.e,

\strut 

$\ \ \ \ \ \ \ q_{2k-1}^{2k+1}=(2N-1)+p_{2k-2}^{2k},$ \ 

$\ \ \ \ \ \ \ q_{2k-3}^{2k+1}=(2N-1)p_{2k-2}^{2k}+p_{2k-4}^{2k},$

...$...,$ \ 

$\ \ \ \ \ \ \ q_{3}^{2k+1}=(2N-1)p_{4}^{2k}+p_{2}^{2k}$ 

and 

$\ \ \ \ \ \ \ q_{1}^{2k+1}=(2N-1)p_{2}^{2k}+p_{0}^{2k}.$ \ \ \ 

$\square $
\end{proposition}

\strut

\strut 

\begin{example}
Suppose that $N=2.$ and let $p_{0}^{2}=4$ and $q_{1}^{3}=3+p_{0}^{2}=3+4=7.$
Put

\strut 

\begin{center}
$G^{8}=X_{8}+p_{6}^{8}X_{6}+p_{4}^{8}X_{4}+p_{2}^{8}X_{4}+p_{0}^{8}e.$
\end{center}

\strut 

Then, by the previous proposition, we have that

\strut 

\begin{center}
$p_{6}^{8}=3+q_{5}^{7},$ \ \ $p_{4}^{8}=3q_{5}^{7}+q_{3}^{7},$ \ \ $%
p_{2}^{8}=3q_{3}^{7}+q_{1}^{7}$ \ \ and \ $p_{0}^{8}=4q_{1}^{7}.$
\end{center}

\strut 

Similarly, by the previous proposition,

\strut 

\ \ \ \ $\ q_{5}^{7}=3+p_{4}^{6},$ \ \ \ $q_{3}^{7}=3p_{4}^{6}+p_{2}^{6}$ \
\ \ \ and \ \ $\ q_{1}^{7}=3p_{2}^{6}+p_{0}^{6},$

$\ \ \ \ \ \ p_{4}^{6}=3+q_{3}^{5},$ \ \ \ $p_{2}^{6}=3q_{3}^{5}+q_{1}^{5}$
\ \ \ \ \ and \ \ \ \ $p_{0}^{6}=4q_{1}^{5},$

\strut 

$\ \ \ \ \ \ q_{3}^{5}=3+p_{2}^{4}$ \ \ \ \ \ \ \ and \ \ \ \ \ \ $%
q_{1}^{5}=3p_{2}^{4}+p_{2}^{4},$

\strut 

$\ \ \ \ \ \ p_{2}^{4}=3+q_{1}^{3}$ \ \ \ \ \ \ \ \ \ and \ \ \ \ \ \ \ \ \ $%
p_{0}^{4}=4q_{1}^{3},$

\strut 

and

\strut 

$\ \ \ \ \ \ \ q_{1}^{3}=3+p_{0}^{2}=7.$

\strut 

Therefore, combining all information,

\strut 

\begin{center}
$G^{8}=X_{8}+22\,X_{6}+202\,X_{4}+744\,X_{2}+1316\,e.$
\end{center}
\end{example}

\strut

We have the following diagram with arrows which mean that

\begin{center}
$\swarrow \swarrow $ \ : \ $(2N-1)+[$former term$]$

$\searrow $ \ \ \ \ : \ \ $(2N-1)\cdot [$former term$]$

$\swarrow $ \ \ \ \ : \ \ $\cdot +[$former term$]$
\end{center}

and

\begin{center}
$\searrow \searrow $ \ : \ \ $(2N)\cdot [$former term$].$
\end{center}

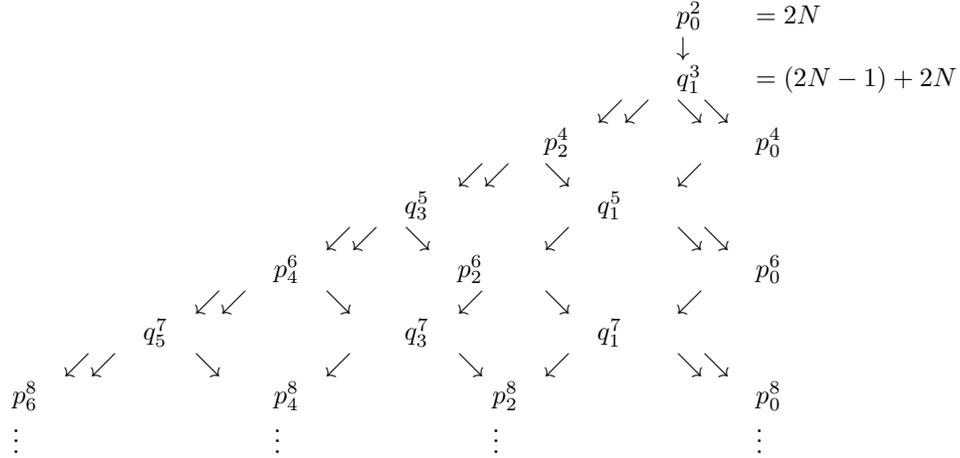
\strut

\begin{center}
$
\begin{array}{llllllllllll}
&  &  &  &  &  &  &  &  &  & p_{0}^{2} & =2N \\ 
&  &  &  &  &  &  &  &  &  & \downarrow  &  \\ 
&  &  &  &  &  &  &  &  &  & q_{1}^{3} & =(2N-1)+2N \\ 
&  &  &  &  &  &  &  &  & \swarrow \swarrow  & \searrow \searrow  &  \\ 
&  &  &  &  &  &  &  & p_{2}^{4} &  &  & p_{0}^{4} \\ 
&  &  &  &  &  &  & \swarrow \swarrow  & \searrow  &  & \swarrow  &  \\ 
&  &  &  &  &  & q_{3}^{5} &  &  & q_{1}^{5} &  &  \\ 
&  &  &  &  & \swarrow \swarrow  & \searrow  &  & \swarrow  &  & \searrow
\searrow  &  \\ 
&  &  &  & p_{4}^{6} &  &  & p_{2}^{6} &  &  &  & p_{0}^{6} \\ 
&  &  & \swarrow \swarrow  &  & \searrow  &  & \swarrow  & \searrow  &  & 
\swarrow  &  \\ 
&  & q_{5}^{7} &  &  &  & q_{3}^{7} &  &  & q_{1}^{7} &  &  \\ 
& \swarrow \swarrow  &  & \searrow  &  & \swarrow  &  & \searrow  & \swarrow 
&  & \searrow \searrow  &  \\ 
p_{6}^{8} &  &  &  & p_{4}^{8} &  &  & \text{ \ \ \ }p_{2}^{8} &  &  &  & 
p_{0}^{8} \\ 
\vdots  &  &  &  & \vdots  &  &  & \text{ \ \ \ }\vdots  &  &  &  & \vdots 
\end{array}
$
\end{center}

\strut 

Recall that Nica and Speicher defined the even random variable in a $*$%
-probability space. Let $(B,\tau _{0})$ be a $*$-probability space, where $%
\tau _{0}:B\rightarrow \Bbb{C}$ is a linear functional satisfying that $\tau
_{0}\left( b^{*}\right) =\overline{\tau _{0}(b)},$ for all $b\in B,$ and let 
$b\in (B,\tau _{0})$ be a random variable. We say that the random variable $%
b\in (B,\tau _{0})$ is even if it is self-adjoint and it satisfies the
following moment relation ;

\strut 

\begin{center}
$\tau _{0}\left( b^{n}\right) =0,$ whenever $n$ is odd.
\end{center}

\strut 

In [12], we observed the amalgamated evenness and we showed that $b\in
(B,\tau _{0})$ is even if and only if

\strut 

\begin{center}
$k_{n}\left( b,...,b\right) =0,$ whenever $n$ is odd.
\end{center}

\strut \strut 

By the previous observation, we can get that

\strut 

\begin{theorem}
Let $G\in \left( A,\tau \right) $ be the generating operator. Then the
moment series of $G$ is

\strut 

\begin{center}
$\tau \left( G^{n}\right) =\left\{ 
\begin{array}{lll}
0 &  & \text{if }n\text{ is odd} \\ 
&  &  \\ 
p_{0}^{n} &  & \text{if }n\text{ is even,}
\end{array}
\right. $
\end{center}

\strut 

for all $n\in \Bbb{N}.$
\end{theorem}

\strut

\begin{proof}
Assume that $n$ is odd. Then 

\strut 

\begin{center}
$G^{n}=X_{n}+q_{n-2}^{n}X_{n-2}+...+q_{3}^{n}X_{3}+q_{1}^{n}X_{1}.$
\end{center}

\strut 

So, $G^{n}$ does not contain the $e$-terms. Therefore,

\strut 

\begin{center}
$\tau \left( G^{n}\right) =\tau \left(
X_{n}+q_{n-2}^{n}X_{n-2}+...+q_{3}^{n}X_{3}+q_{1}^{n}X_{1}\right) =0.$
\end{center}

\strut 

Assume that $n$ is even. Then

\strut 

\begin{center}
$G^{n}=X_{n}+p_{n-2}^{n}X_{n-2}+...+p_{2}^{n}X_{2}+p_{0}^{n}e.$
\end{center}

\strut 

So, we have that

\strut 

\begin{center}
$\tau (G^{n})=\tau \left(
X_{n}+p_{n-2}^{n}X_{n-2}+...+p_{2}^{n}X_{2}+p_{0}^{n}e\right) =p_{0}^{n}.$
\end{center}
\end{proof}

\strut 

Remark that the $n$-th moments of the generating operator in $(A,\tau )$ is
totally depending on the recurrence relation.\strut 

\strut 

\begin{corollary}
Let $G\in (A,\tau )$ be the generating operator. Then $G$ is even. $\square $
\end{corollary}

\strut 

\begin{corollary}
Let $G\in (A,\tau )$ be the generating operator. Then

\strut 

\begin{center}
$M_{G}(z)=\sum_{n=1}^{\infty }p_{0}^{2n}\,z^{2n}\in \Theta _{1}.$
\end{center}
\end{corollary}

\strut 

\begin{proof}
Since all odd moments of $G$ vanish, $coef_{n}\left( M_{G}\right) =0,$ for
all odd integer $n.$ By the previous theorem, we can get the result.
\end{proof}

\strut 

\strut 

\strut 

\section{The R-transform of the Generating Operator of $L(F_{N})$}

\strut 

\strut 

In this chapter, we will compute the R-transform of the generating operator $%
G$ in $(A,\tau )\equiv \left( L(F_{N}),\tau \right) .$ We can get the
R-transform by using the M\"{o}bius inversion, which we considered in
Chapter 1. Notice that, by the evenness of $G,$ we have that all od
cumulants of $G$ vanish. i.e,

\strut 

\begin{center}
$k_{n}\left( \underset{n-times}{\underbrace{G,........,G}}\right) =0,$
whenever $n$ is odd.
\end{center}

\strut 

This shows that the nonvanishing $n$-th coefficients of the R-transform of $%
G,$ $R_{G}$, are all even coefficients.

\strut \strut 

\begin{theorem}
Let $G\in (A,\tau )$ be the generating operator. Then $R_{G}(z)=\sum_{n=1}^{%
\infty }\alpha _{2n}\,\,z^{2n},$ in $\Theta _{1},$ with

\strut 

\begin{center}
$\alpha _{2n}=\underset{l_{1},...,l_{p}\in 2\Bbb{N},\,l_{1}+...+l_{p}=2n}{%
\sum }\,\underset{\pi \in NC_{l_{1},...,l_{p}}(2n)}{\sum }\left(
p_{0}^{l_{1}}...p_{0}^{l_{p}}\right) \mu (\pi ,1_{2n}),$
\end{center}

\strut 

for all $n\in \Bbb{N},$ where

\strut 

\begin{center}
$NC_{l_{1},...,l_{p}}(2n)=\{\pi \in NC(2n):V\in \pi \Leftrightarrow \left|
V\right| =l_{j},\,j=1,...,p\}.$
\end{center}
\end{theorem}

\strut 

\begin{proof}
By the evenness of $G,$ all odd cumulants vanish (See [12]). Fix $n\in \Bbb{N%
},$ an even number. Then

\strut 

$\ \ \ k_{n}\left( \underset{n-times}{\underbrace{G,........G}}\right) =%
\underset{\pi \in NC(n)}{\sum }\tau _{\pi }\left( G,...,G\right) \,\mu (\pi
,1_{n})$

\strut 

$\ \ \ \ \ \ \ \ \ \ \ \ \ \ \ \ =\underset{\pi \in NC^{(even)}(n)}{\sum }%
\tau _{\pi }\left( G,...,G\right) \,\mu (\pi ,1_{n})$

\strut 

where $NC^{(even)}(n)=\{\pi \in NC(n):\pi $ does not contain odd blocks$\},$
by [12]

\strut 

$\ \ \ \ \ \ \ \ \ \ \ \ \ \ \ \ =\underset{\pi \in NC^{(even)}(n)}{\sum }%
\left( \underset{V\in \pi }{\prod }\tau (G^{\left| V\right| })\right) \mu
(\pi ,1_{n})$

\strut 

\strut (2.1)

$\ \ \ \ \ \ \ \ \ \ \ \ \ \ \ \ =\underset{\pi \in NC^{(even)}(n)}{\sum }%
\left( \underset{V\in \pi }{\prod }\,p_{0}^{\left| V\right| }\right) \,\mu
(\pi ,1_{n}).$

\strut 

By [14], we have that

\strut 

\begin{center}
$NC^{(even)}(n)=\underset{l_{1},...,l_{p}\in 2\Bbb{N},\,l_{1}+...+l_{p}=n}{%
\bigsqcup }NC_{l_{1},...,l_{p}}(n)$
\end{center}

\strut 

where $\bigsqcup $ is the disjoint union and

\strut 

\begin{center}
$NC_{l_{1},...,l_{p}}\left( n\right) =\{\pi \in NC^{(even)}(n):V\in \pi
\Leftrightarrow \left| V\right| =l_{j},\,j=1,..,p\}.$
\end{center}

\strut 

(For instance, $NC^{(even)}(6)=NC_{2,2,2}(6)\cup NC_{2,4}(6)\cup NC_{6}(6).$)

\strut 

So, the formular (2.1) goes to

\strut 

$\underset{l_{1},...,l_{p}\in 2\Bbb{N},\,l_{1}+...+l_{p}=n}{\sum }\,%
\underset{\pi \in NC_{l_{1},...,l_{p}}(n)}{\sum }\left(
p_{0}^{l_{1}}...p_{0}^{l_{p}}\right) \mu (\pi ,1_{n}).$
\end{proof}

\strut 

\begin{example}
Let $N=2.$ Then $A=L(F_{2})$ and $G=a+b+a^{-1}+b^{-1},$ where $F_{2}=\,<a,b>.
$ Then

\strut 

$\ \ 
\begin{array}{ll}
k_{4}\left( G,G,G,G\right)  & =\underset{\pi \in NC_{2,2}(4)}{\sum }\left(
p_{0}^{2}p_{0}^{2}\right) \,\mu (\pi ,1_{4})+p_{0}^{4} \\ 
& =-2\left( p_{0}^{2}\right) ^{2}+p_{0}^{4}=-32+28 \\ 
& =-4.
\end{array}
$

\strut and

$\ 
\begin{array}{ll}
k_{6}\left( G,...,G\right)  & =\underset{\pi \in NC_{2,2,2}(6)}{\sum }\left(
p_{0}^{2}p_{0}^{2}p_{0}^{2}\right) \mu (\pi ,1_{6}) \\ 
& \,\,\,\,\,\,\,\,\,\,\,\,\,\,\,\,\,\,\,\,\,\,\,\,\,\,\,\,\,\,\,\,\,\,\,\,\,%
\,\,\,\,\,\,\,\,\,\,\,\,\,+\underset{\pi \in NC_{2,4}(6)}{\sum }\left(
p_{0}^{2}p_{0}^{4}\right) \mu (\pi ,1_{6})+p_{0}^{6} \\ 
& =\left(
2(p_{0}^{2})^{3}+2(p_{0}^{2})^{3}+(p_{0}^{2})^{3}+(p_{0}^{2})^{3}+(p_{0}^{2})^{3}\right) 
\\ 
& \,\,\,\,\,\,\,\,\,\,\,\,\,\,\,\,\,\,\,\,\,\,\,\,\,\,\,\,\,\,\,\,\,\,\,\,\,%
\,\,\,\,\,\,\,\,\,\,\,\,\,+\left| NC_{2,4}(6)\right| \left(
p_{0}^{2}p_{0}^{4}(-1)\right) +p_{0}^{6} \\ 
& =448-672+232=8.
\end{array}
$
\end{example}

\strut 

\strut 

\strut

\strut \textbf{References}

\strut

\label{REF}

\strut

\begin{enumerate}
\item[{[1]}]  R. Speicher, Combinatorial Theory of the Free Product with
Amalgamation and Operator-Valued Free Probability Theory, AMS Mem, Vol 132 ,
Num 627 , (1998).

\item[{\lbrack 2]}]  A. Nica, R-transform in Free Probability, IHP course
note, available at

\texttt{www.math.uwaterloo.ca/\symbol{126}anica.}\strut 

\item[{[3]}]  R. Speicher, Combinatorics of Free Probability Theory IHP
course note,

available at \texttt{www.mast.queensu.ca/\symbol{126}speicher.}\strut 

\item[{[4]}]  A. Nica, D. Shlyakhtenko and R. Speicher, R-cyclic Families of
Matrices in Free Probability, J. of Funct Anal, 188 (2002), 227-271.\strut 

\item[{[5]}]  A. Nica and R. Speicher, R-diagonal Pair-A Common Approach to
Haar Unitaries and Circular Elements, (1995), \texttt{www.mast.queensu.ca/%
\symbol{126}speicher.}\strut 

\item[{[6]}]  D. Shlyakhtenko, Some Applications of Freeness with
Amalgamation, J. Reine Angew. Math, 500 (1998), 191-212.\strut 

\item[{[7]}]  A. Nica, D. Shlyakhtenko and R. Speicher, R-diagonal Elements
and Freeness with Amalgamation, Canad. J. Math. Vol 53, Num 2, (2001)
355-381.\strut 

\item[{[8]}]  A. Nica, R-transforms of Free Joint Distributions and
Non-crossing Partitions, J. of Func. Anal, 135 (1996), 271-296.\strut 

\item[{[9]}]  D.Voiculescu, K. Dykemma and A. Nica, Free Random Variables,
CRM Monograph Series Vol 1 (1992).\strut 

\item[{[10]}]  D. Voiculescu, Operations on Certain Non-commuting
Operator-Valued Random Variables, Ast\'{e}risque, 232 (1995), 243-275.\strut 

\item[{[11]}]  D. Shlyakhtenko, A-Valued Semicircular Systems, J. of Funct
Anal, 166 (1999), 1-47.\strut 

\item[{[12]}]  I. Cho, Amalgamated Boxed Convolution and Amalgamated
R-transform Theory (preprint).\strut 

\item[{\lbrack 13]}]  I. Cho, Compatibility of a noncommutative probability
space and an amalgamated noncommutative probability space, preprint

\item[{\lbrack 14]}]  I. Cho, An Example of Moment Series under the
Compatibility, Preprint

\item[{\lbrack 15]}]  F. Radulescu, Singularity of the Radial Subalgebra of $%
L(F_{N})$ and the Puk\'{a}nszky Invariant, Pacific J. of Math, vol. 151, No
2 (1991)\strut , 297-306.\strut \strut 
\end{enumerate}

\strut

\end{document}